\documentclass{amsart}
\numberwithin{equation}{section}

\newtheorem{thm}{Theorem}[section]
\newtheorem{lem}[thm]{Lemma}

\newtheorem{defin}[thm]{Definition}
\newtheorem{rem}[thm]{Remark}

\begin{document}
\title{Initial-boundary value problem}

\author{Oqila Muhiddinova}
\address{Institute of Mathematics, Uzbekistan Academy of Science}
\curraddr{Institute of Mathematics, Uzbekistan Academy of Science,
Tashkent, 81 Mirzo Ulugbek str. 100170} \email{oqila1992@mail.ru}

\small

\title[Initial-boundary value problem]
{Initial-boundary value problem for a subdiffusion equation with
the Caputo derivative}

\begin{abstract}

We investigate an initial-boundary value problem for a
time-fractional subdiffusion equation with the Caputo derivatives
on $N$-dimensional torus by the classical Fourier method. Since
our solution is established on the eigenfunction expansion of
elliptic operator, the method proposed in this article can be used
to an arbitrary domain and an elliptic operator with variable
coefficients. It should be noted that the conditions for the
existence of a solution to the initial-boundary value problem
found in the article cannot be weakened, and the article provides
a corresponding example.

\vskip 0.3cm \noindent {\it AMS 2000 Mathematics Subject
Classifications} :
Primary 35R11; Secondary 74S25.\\
{\it Key words}: Time-fractional subdiffusion equation, the Caputo
derivatives, initial-boundary value problem, Fourier method,
Liouville spaces.
\end{abstract}+-

\maketitle

\section{Main result}

The fractional derivative in the sense of Caputo of order $0<\rho
<1$ of the function $h (t)$ defined on $[0, \infty)$ has the
form (see, for example, \cite{PSK}, p. 14)
$$
D_t^\rho h(t)= \frac{1}{\Gamma(1-\rho)} \int_0^t \frac{\frac{d
}{d\tau}h(\tau)d\tau}{(t-\tau)^{\rho}}, \quad t>0,
$$
provided the right-hand side exists. Here $\Gamma(\sigma)$ is
Euler's gamma function. If in this definition we interchange
differentiation and fractional integration, then we get the
definition of the Riemann-Liouville derivative:
$$
\partial_t^\rho h(t)= \frac{1}{\Gamma(1-\beta)}\frac{d}{dt}
\int_0^t \frac{h(\tau)d\tau}{(t-\tau)^{\rho}}, \quad t>0.
$$

Note that if $\rho=1$, then fractional derivatives coincides with
the ordinary classical derivative of the first order: 
$\partial_th(t) = D_t h(t)= \frac{d}{dt} h(t)$.

Let $\mathbb{T}^N$ be $N$ -dimensional torus: $\mathbb{T}^N=(-\pi,
\pi]^N$, $N\geq 1$. We define by $C(\mathbb{T}^N)$ and
$C^2(\mathbb{T}^N)$ a class of $2\pi$-periodic on each variable
$x_j$ functions $v(x)$ from $C(\mathbb{T}^N)$ and
$C^2(\mathbb{T}^N)$ correspondingly. Let $A$ stand for a positive
operator, defined on $C^2(\mathbb{T}^N)$ and acting as $A
v(x)=-\Delta v(x)$, where $\Delta$ is the Laplace operator.

Let $\rho\in(0,1) $ be a constant number. Consider the
initial-boundary value problem
\begin{equation}\label{eq}
D_t^\rho u(x,t) + Au(x,t) = f(x,t), \quad x\in \mathbb{T}^N, \quad
0<t\leq T,
\end{equation}

\begin{equation}\label{in}
u(x,0) = \varphi(x), \quad x\in \mathbb{T}^N,
\end{equation}
where $f$ and $\varphi$ are given continuous functions.

\begin{defin} A function $u(x,t)\in C(\mathbb{T}^N\times
[0,T])$ with the properties $D_t^\rho u(x,t), A(x,D)u(x,t)\in
C(\mathbb{T}^N\times (0,T])$, and satisfying the conditions of
problem (\ref{eq}) - (\ref{in}) is called \textbf{the classical
solution} of the initial-boundary value problem.
\end{defin}

Initial-boundary value problem (\ref{eq}) - (\ref{in}) for various
elliptic operators $A$ have been considered by a number of authors
using different methods (see, for example, handbook \cite{Koch}).
In the book of A.A. Kilbas et al. \cite{Kil} (Chapter 6) there is
a survey of works published before 2006. The case of one spatial
variable $x\in \mathbb{R}$ and subdiffusion equation with the
elliptical part $u_{xx}$ were considered  for example in the
monograph of A. V. Pskhu \cite{PSK} (Chapter 4, see references
thesein).

The main method used in this work is the Fourier method. As far as
we know, in all previous papers this method was used to prove the
existence of a generalized solution of initial-boundary value
problems for subdiffusion equations. So for example, in the paper
of Yu. Luchko \cite{Luch} the author constructed solutions by the
eigenfunction expansion in the case of $f = 0$ and discussed the
unique existence of the generalized solution to problem (\ref{eq})
- (\ref{in}) with the Caputo derivative.  In an arbitrary
$N$-dimensional domain $\Omega$ initial-boundary value problems
for subdiffusion equations (the fractional part of the equation is
a multi-term and initial conditions are non-local) with the Caputo
derivatives has been investigated by M. Ruzhansky et al.
\cite{Ruz}. The authors proved the existence and uniqueness of the
generalized solution to the problem by the Fourier method. The
authors of the paper \cite{Yama11} investigated both problem
(\ref{eq}) - (\ref{in}) and the corresponding backward problem
with an arbitrary elliptic operator of the second order.  To prove
the existence and uniqueness of the generalized solution the
authors applied the Fourier method.

The authors of papers \cite{AO} and \cite{AO1} used the Fourier
method to construct a classical solution of the subdiffusion
equations with the Riemann-Liouville derivative and various
elliptic operators.

Let us introduce some concepts and formulate the main result of
the work.

Let $\hat{A}$ stand for the closure  of operator $A$ in
$L_2(\mathbb{T}^N)$. Then $\hat{A}$ is selfadjoint and it has a
complete (in $L_2(\mathbb{T}^N)$) set of eigenfunctions $\{\gamma
e^{inx}\}$, $\gamma = \gamma(N)=(2\pi)^{-N/2}$, $n\in
\mathbb{Z}^N$ and corresponding eigenvalues
$|n|^2=n_1^2+n_2^2+...+n_N^2$. Therefore, by virtue of J. von
Niemann theorem, for any $\tau> 0$ one can introduce the power of
operator $\hat{A}$ as $\hat{A}^\tau g(x)=\sum\limits_{n\in
\mathbb{Z}^N} |n|^\tau g_n e^{inx}$, where $g_n$ are Fourier
coefficients:
$$
g_n=(2\pi)^{-N}\int\limits_{\mathbb{T}^N} g(x) e^{-inx} dx.
$$
The domain of definition of this operator is defined from the
condition $\hat{A}^\tau g(x)\in L_2(\mathbb{T}^N)$ and has the
form
$$
D(\hat{A}^\tau)=\{g\in L_2(\mathbb{T}^N): \sum\limits_{n\in
\mathbb{Z}^N} |n|^{2\tau} |g_n|^2 < \infty\}.
$$
If we denote by $L_2^a(\mathbb{T}^N)$, $a> 0$, the Liouville space
with the norm
\begin{equation}\label{T}
||g||^2_{L_2^a(\mathbb{T}^N)}=\big|\big|\sum\limits_{n\in
\mathbb{Z}^N}(1+|n|^2)^{\frac{a}{2}}g_n
e^{inx}\big|\big|^2_{L_2(\mathbb{T}^N)}=\sum\limits_{n\in
\mathbb{Z}^N}(1+|n|^2)^a|g_n|^2,
\end{equation}
then it is not hard to see that $ D(\hat{A}^\tau)= L_2^{\tau
m}(\mathbb{T}^N)$.

Let $E_{\rho, \mu}$ be  the two-parametric Mittag-Leffler
function:
$$
E_{\rho, \mu}(t)= \sum\limits_{k=0}^\infty \frac{t^k}{\Gamma(\rho
k+\mu)}.
$$

Here is the main result.

\begin{thm}\label{main} Let $a > \frac{N}{2}$ and $\varphi\in L^{a}_2(\mathbb{T}^N)$. Moreover, let $f(x,t)\in
L^a_2(\mathbb{T}^N)$ for $t\in [0, T]$
 and $||f(\cdot, t)||^2_{L_2^a(\mathbb{T}^N)}\in C[0, T]$. Then there exists a solution of initial-boundary
value problem (\ref{eq}) - (\ref{in}) and it has the form
\begin{equation}\label{solution}
u(x,t)=\sum\limits_{n\in\mathbb{Z}^N} \bigg[\varphi_n E_{\rho, 1}
(-|n|^2 t^\rho)+\int\limits_0^t f_n(t-\xi)\xi^{\rho-1} E_{\rho,
\rho}(-|n|^2\xi^\rho) d\xi\bigg]e^{inx},
\end{equation}
which absolutely and uniformly converges on $x\in \mathbb{T}^N$
and for each $t\in (0, T]$, where $\varphi_n$ and $f_n(t)$ are
corresponding Fourier coefficients. Moreover, the series obtained
after applying term-wise the operators $D_t^\rho$ and $A$ also
converge absolutely and uniformly on $x\in \mathbb{T}^N$ and for
each $t\in (0, T]$.

\end{thm}
\begin{rem}\label{rem}
When $a>\frac{N}{2}$, according to the Sobolev embedding theorem,
all functions in $L^{a}_2(\mathbb{T}^N)$ are $2\pi$-periodic
continuous functions. The fulfillment of the inverse inequality $a
\leq \frac{N}{2}$, admits the existence of unbounded functions in
$L^{a}_2(\mathbb{T}^N)$ (see, for example, \cite{AAP}). Therefore,
condition  $a>\frac{N}{2}$ for function $f$ of this theorem is not
only sufficient for the statement to be hold, but it is also
necessary.
\end{rem}

A result similar to Theorem \ref{main} was obtained in the recent
paper \cite{AO1} for a subdiffusion equation with the
Riemann-Liouville derivative. But the condition found for the
initial function $\varphi(x)$ in that work is less restrictive.
For example, in the one-dimensional case, it has been proved that
it is sufficient to require $\varphi(x)\in C(\mathbb{T}^1)$. As
will be proved at the end of this paper, even Holder continuous
with exponent $a=\frac{1}{2}$ for the initial function is not
sufficient for the validity of Theorem \ref{main} in the
one-dimensional case.

In conclusion, note that since our solution is established on the
eigenfunction expansion of elliptic operator, the method proposed
in this article can be used to an arbitrary domain and an elliptic
operator with variable coefficients.

\section{Uniqueness}

In this section we prove the uniqueness of the solution to problem
(\ref{eq})-(\ref{in}). Note that in order for problem
(\ref{eq})-(\ref{in}) to have a unique solution, it is sufficient
that the functions $f$ and $\varphi$ be continuous.

Suppose that initial-boundary value problem (\ref{eq}) -
(\ref{in}) has two classical solutions $u_1(x,t)$ and $u_2(x,t)$.
Our aim is to prove that $u(x,t)=u_1(x,t)-u_2(x,t)\equiv 0$. Since
the problem is linear, then we have the following homogenous
problem for $u(x,t)$:
\begin{equation}\label{eq1}
D_t^\rho u(x,t) + A\,u(x,t) = 0, \quad x\in \mathbb{T}^N, \quad
t>0;
\end{equation}
\begin{equation}\label{in1}
u(x,0) = 0, \quad x\in \mathbb{T}^N.
\end{equation}

Let $u(x,t)$ be a solution of problem (\ref{eq1})-(\ref{in1}).
Consider the function
\begin{equation}\label{w}
w_n(t)=\int\limits_{\mathbb{T}^N} u(x,t)e^{inx}dx, \quad n\in
\mathbb{Z}^N.
\end{equation}
By virtue of equation (\ref{eq1}), we can write
$$
D_t^\rho w_n(t)=\int\limits_{\mathbb{T}^N} D_t^\rho
u(x,t)e^{inx}dx= -\int\limits_{\mathbb{T}^N} A \, u(x,t)e^{inx}dx,
\quad t>0,
$$
or, integrating by parts,
$$
D_t^\rho w_k(t)=-\int\limits_{\mathbb{T}^N}
u(x,t)A\,e^{inx}dx=-|n|^2 \int\limits_{\mathbb{T}^N}
u(x,t)e^{inx}dx = -|n|^2 w_n(t), \quad t>0.
$$
Using in (\ref{w}) the homogenous initial condition (\ref{in1}),
we have the following Cauchy problem for $w_n(t)$:
$$
D_t^\rho w_n(t) +|n|^2 w_n(t)=0,\quad t>0; \quad w_n(0)=0.
$$
This problem has the unique solution; therefore, the function
defined by (\ref{w}), is identically zero: $w_n(t)\equiv 0$ (see
for example, \cite{PSK} p. 17 , \cite{ACT}). From completeness in
$L_2(\mathbb{T}^N)$ of the system of eigenfunctions $\{e^{inx}\}$,
we have $u(x,t) = 0$ for all $x\in \mathbb{T}^N$ and $t>0$. Hence
the uniqueness is proved.

\section{Existence}

Proof of existence based on the following lemma (see M.A.
Krasnoselski et al. \cite{Kra}, p. 453), which is a simple
corollary of the Sobolev embedding theorem.

\begin{lem}\label{CL} Let $\sigma > 1+\frac{N}{4}$. Then for any $|\alpha|\leq 2$
operator $D^\alpha (\hat{A}+1)^{-\sigma}$ (completely)
continuously maps from $L_2(\mathbb{T}^N)$ into $C(\mathbb{T}^N)$
and moreover the following estimate holds true
\begin{equation}\label{CL1}
||D^\alpha (\hat{A}+1)^{-\sigma} g||_{C(\mathbb{T}^N)} \leq C
||g||_{L_2(\mathbb{T}^N)}.
\end{equation}

\end{lem}
Proof of this lemma see in \cite{AO1}.

One can easily verify that the function (\ref{solution}) formally
satisfies the conditions of problem (\ref{eq})-(\ref{in}). In
order to prove that function (\ref{solution}) is actually a
solution to the problem, it remains to substantiate this formal
statement, i.e. show that the operators $A$ and $D_t^\rho$ can be
applied term by term to the series (\ref{solution}). To do this we
remind the following asymptotic estimate of the Mittag-Leffler
function with a sufficiently large negative argument $t$ and an
arbitrary complex number $\mu$ (see, for example, \cite{Dzh66}, p.
134)
\begin{equation}\label{m1}
|E_{\rho, \mu}(-t)|\leq \frac{C}{1+ t}, \quad t>0.
\end{equation}
We will also use a coarser estimate with a positive $\lambda$ and
$0<\varepsilon<1$:
\begin{equation}\label{m2}
|t^{\rho-1} E_{\rho,\mu}(-\lambda t^\rho)|\leq
\frac{Ct^{\rho-1}}{1+\lambda t^\rho}\leq C \lambda^{\varepsilon-1}
t^{\varepsilon\rho-1}, \quad t>0,
\end{equation}
which is easy to verify. Indeed, let $t^\rho\lambda<1$, then $t<
\lambda^{-1/\rho}$ and
$$
t^{\rho -1} = t^{\rho-\varepsilon\rho} t^{\varepsilon\rho-1} <
\lambda^{\varepsilon-1}t^{\varepsilon\rho-1}.
$$
If $t^\rho\lambda\geq 1$, then $\lambda^{-1}\leq t^\rho$ and
$$
\lambda^{-1} t^{-1}=\lambda^{-1+\varepsilon}
\lambda^{-\varepsilon} t^{-1}\leq
\lambda^{\varepsilon-1}t^{\varepsilon\rho-1}.
$$

Note the series (\ref{solution}) is in fact the sum of two series.
Consider the following partial sums of the first series:
\begin{equation}\label{S1}
S^1_k(x,t)=\sum\limits_{|n|^2<k}  E_{\rho, 1}(-|n|^2
t^\rho)\varphi_n \, e^{inx},
\end{equation}
and suppose that function $\varphi$ satisfies the condition of
Theorem \ref{main}, i.e. for some $\tau> \frac{N}{4}$
$$
\sum\limits_{n\in\mathbb{Z}^N} |n|^{4\tau} |\varphi_n|^2 \leq
C_\varphi<\infty.
$$
Since $\hat{A}^{-\tau-1} e^{inx} = |n|^{-2(\tau+1)} e^{inx}$, we
may rewrite the sum (\ref{S1}) as
$$
S^1_k(x,t)=\hat{A}^{-\tau-1}\sum\limits_{|n|^2<k}
E_{\rho,1}(-|n|^2 t^\rho)\varphi_n \, |n|^{2(\tau+1)}\, e^{inx}.
$$
Therefore by virtue of Lemma \ref{CL} one has
$$
||D^\alpha S^1_k||_{C(\mathbb{T}^N)}=||D^\alpha
\hat{A}^{-\tau-1}\sum\limits_{|n|^2<k} E_{\rho,1}(-|n|^2
t^\rho)\varphi_n \, |n|^{2(\tau+1)}\,
e^{inx}||_{C(\mathbb{T}^N)}\leq
$$
\begin{equation}\label{S11}
\leq C ||\sum\limits_{|n|^2<k}  E_{\rho,1}(-|n|^2 t^\rho)\varphi_n
\, |n|^{2(\tau+1)}\, e^{inx}||_{L_2(\mathbb{T}^N)}.
\end{equation}
Using the orthonormality of the system $\{e^{inx}\}$, we will have
\begin{equation}\label{S2}
||D^\alpha S^1_k||^2_{C(\mathbb{T}^N)}\leq C \sum\limits_{|n|^2<k}
\big| E_{\rho,1}(-|n|^2 t^\rho)\varphi_n \,
|n|^{2(\tau+1)}\big|^2.
\end{equation}
Application of estimate (\ref{m1}) and inequality $|n|^2 (1+ |n|^2
t^\rho)^{-1}<t^{-\rho}$ gives
$$
\sum\limits_{|n|^2<k} \big|E_{\rho,1}(-|n|^2 t^\rho)\varphi_n \,
|n|^{2(\tau+1)}\big|^2\leq C t^{-2\rho}\sum\limits_{|n|^2<k}
|n|^{4\tau}|\varphi_n|^2\leq C t^{-2\rho}C_\varphi.
$$
Therefore we can rewrite the estimate (\ref{S2}) as
$$
||D^\alpha S^1_k||^2_{C(\mathbb{T}^N)}\leq C t^{-2\rho} C_\varphi.
$$

This implies uniformly on $x\in\mathbb{T}^N$ convergence of the
differentiated sum (\ref{S1}) with respect to the variables $x_j$
for each $t\in (0,T]$. On the other hand, the sum (\ref{S11})
converges for any permutation of its members as well, since these
terms are mutually orthogonal. This implies the absolute
convergence of the differentiated sum (\ref{S1}) on the same
interval $t\in (0,T]$.

Now we consider the second part of the series (\ref{solution}):
\begin{equation}\label{S2}
S^2_k(x,t)=\sum\limits_{|n|^2<k} \int\limits_0^t
f_n(t-\xi)\xi^{\rho-1} E_{\rho,\rho}(-|n|^2 \xi^\rho)\,d\xi\,
e^{inx}
\end{equation}
and suppose that function  $f(x,t)$  satisfies all the conditions
of Theorem \ref{main}, i.e. the following series converges
uniformly on $t\in [0, T]$ for some $\tau> \frac{N}{4}$:
$$
\sum\limits_{n\in\mathbb{Z}^N} |n|^{4\tau} |f_n(t)|^2 \leq
C_f<\infty.
$$
We choose a small $\varepsilon>0$ in such a way, that
$\tau+1-\varepsilon> 1+\frac{N}{4}$.  Since
$\hat{A}^{-\tau-1+\varepsilon} e^{inx} =
|n|^{-2(\tau+1-\varepsilon)} e^{inx}$, we may rewrite the sum
(\ref{S2}) as
$$
S^2_k(x,t)=\hat{A}^{-\tau-1+\varepsilon}\sum\limits_{|n|^2<k}
\int\limits_0^t f_n(t-\xi)\xi^{\rho-1} E_{\rho,\rho}(-|n|^2
\xi^\rho)\,d\xi\,|n|^{2(\tau+1-\varepsilon)} e^{inx}.
$$
Then by virtue of Lemma \ref{CL} one has
$$
||D^\alpha S^2_k||_{C(\mathbb{T}^N)}=||D^\alpha
\hat{A}^{-\tau-1+\varepsilon}\sum\limits_{|n|^2<k} \int\limits_0^t
f_n(t-\xi)\xi^{\rho-1} E_{\rho,\rho}(-|n|^2
\xi^\rho)\,d\xi\,|n|^{2(\tau+1-\varepsilon)}
e^{inx}||_{C(\mathbb{T}^N)}\leq
$$
\begin{equation}\label{S4}
\leq C \big|\big|\sum\limits_{|n|^2<k} \int\limits_0^t
f_n(t-\xi)\xi^{\rho-1} E_{\rho,\rho}(-|n|^2
\xi^\rho)\,d\xi\,|n|^{2(\tau+1-\varepsilon)}
e^{inx}\big|\big|_{L_2(\mathbb{T}^N)}.
\end{equation}
Using the orthonormality of the system $\{e^{inx}\}$, we will have
$$
||D^\alpha S^2_k||^2_{C(\mathbb{T}^N)}\leq C \sum\limits_{|n|^2<k}
\big|\int\limits_0^t f_n(t-\xi)\xi^{\rho-1} E_{\rho,\rho}(-|n|^2
\xi^\rho)\,d\xi\,|n|^{2(\tau+1-\varepsilon)}\big|^2.
$$
Now we use estimate (\ref{m2}) and apply the generalized Minkowski
inequality. Then
$$
||D^\alpha S^2_k||^2_{C(\mathbb{T}^N)}\leq C\bigg(\int\limits_0^t
\xi^{\varepsilon\rho-1} \big(\sum\limits_{|n|^2<k}
|n|^{4\tau}|f_n(t-\xi)|^2\big)^{1/2} d\xi \bigg)^2\leq C\cdot C_f,
$$
where $C$ depends on $T$ and $\varepsilon$. Hence, using the same
argument as above, we see that the differentiated sum (\ref{S2})
with respect to the variables $x_j$ converges absolutely and
uniformly on $(x,t)\in \mathbb{T}^N\times [0,T]$.

Further, from equation (\ref{eq}) one has
$$
D_t^\rho \big(S^1_k+S^2_k\big)=  -A
\big(S^1_k+S^2_k\big)+\sum\limits_{|n|^2<k} f_n(t)e^{inx}.
$$
Absolutely and uniformly convergence of the latter  series can be
proved as above.

Thus Theorem \ref{main} is completely proved.

\section{Counterexample}

In this section, we will discuss the importance of the condition
$a>\frac{N}{2}$ of Theorem \ref{main}. When $a>\frac{N}{2}$ all
functions in $L_2^a(\mathbb{T}^N)$ belong to $C(\mathbb{T}^N)$. As
noted in Remark \ref{rem}, if $a=\frac{N}{2}$, then in the class
$L_2^a(\mathbb{T}^N)$ there exist unbounded functions, as a
consequence of which problem (\ref{eq})-(\ref{in}) certainly does
not have classical solutions. Therefore, the question naturally
arises: is it possible to replace, for example, condition
$\varphi\in L_2^a(\mathbb{T}^N)$, $a>\frac{N}{2}$,  of Theorem
\ref{main} by condition
\begin{equation}\label{class}
\varphi\in
L_2^{\frac{N}{2}}(\mathbb{T}^N)\cap C(\mathbb{T}^N)?
\end{equation}
The following example answers this question in the
negative.

Let $N=1$. In the class of periodic functions, we seek a solution
to the following problem
\begin{equation}\label{prob1}
\left\{
\begin{aligned}
&D_t^\rho u(x,t) -u_{xx}(x, t) = 0,\quad 0<t\leq T;\\
&u(x, o) = \varphi(x).
\end{aligned}
\right.
\end{equation}

If $\varphi$ satisfies the condition of Theorem \ref{main}, then
the unique solution of the problem has the form
\begin{equation}\label{solution1}
u(x,t)=\lim\limits_{k\rightarrow \infty}\sum\limits_{|n|\leq k}
\varphi_n E_{\rho, 1} (-|n|^2 t^\rho)\,e^{inx}.
\end{equation}

Recall, that the classes $C^a(\mathbb{T}^1)$ are usually defined
as follows: $2\pi$-periodic function $\varphi\in
C^a(\mathbb{T}^1)$ if and only if
\[
|\varphi(x) - \varphi(y)|\leq C |x-y|^a.
\]
Consider the following function, first studied by Hardy and
Littlewood (see \cite{Zyg}, proof of Theorem (3.10)):
\[
\Phi(x) =\sum\limits_{n=1}^\infty \frac{e^{in\ln n}}{n}e^{inx}.
\]
The real and imaginary parts of this function belong to the class
$C^{\frac{1}{2}}(\mathbb{T}^1)$ (see \cite{Zyg}, Chapter V,
paragraph 4). Set $\varphi(x)= \Re(\Phi(x))$, where $\Re(z)$ is a
real part of the complex number $z$. Then $\varphi\in
C^{\frac{1}{2}}(\mathbb{T}^1)$ and it is not hard to see, that
\[
\sum\limits_{n=1}^\infty
\sqrt{(\varphi^c_n)^2+(\varphi^s_n)^2}=+\infty,
\]
where $\varphi^c_n$ are the coefficients of the function
$\varphi(x)$ in $\cos nx$ and $\varphi^s_n$ - in terms of $\sin
nx$. Obviously, the function $\varphi(x)$ also belongs to the
class (\ref{class}) and if we denote
$\varphi_n=\frac{1}{2}(\varphi^c_n - i \varphi^s_n)$ and
$\varphi_{-n}= \overline{\varphi_n}$, then
\begin{equation}\label{phi}
\sum\limits_{n\in \mathbb{Z}^1} |\varphi_n|=+\infty.
\end{equation}

Suppose that the solution to problem (\ref{prob1}) has the form
(\ref{solution1}). Let us show that the series (\ref{solution1})
differentiated twice with respect to the variable $x$ does not
converge absolutely, that is, for the function $\varphi(x)$,
defined above, the statement of Theorem \ref{main} does not hold.
Indeed, set
\[
(u_{k_0})_{xx}(x,t)=-\sum\limits_{k_0\leq |n|} \varphi_n E_{\rho,
1} (-|n|^2 t^\rho)\,|n|^2\,e^{inx}.
\]
In order for this series to converge uniformly and absolutely with
respect to $x\in\mathbb{T}^1$ and for each $t\in(0,T]$ it is
necessary that the number series
\[
U_{k_0}=\sum\limits_{k_0\leq |n|} |\varphi_n E_{\rho, 1} (-|n|^2
t^\rho)|\,|n|^2
\]
converge for some $t>t_0>0$. Now we remind the following
asymptotic estimate of the Mittag-Leffler function with a
sufficiently large negative argument
\[
E_{\rho, 1} (-|n|^2
t^\rho)=\frac{1}{\Gamma(1-\rho)}\cdot\frac{1}{|n|^2
t^\rho}+O\bigg(\frac{1}{|n|^2 t^\rho}\bigg)^2.
\]
Hence, for sufficiently large $k_0$, we have
\[
U_{k_0}= \frac{1}{t^{\rho}\Gamma(1-\rho)}\cdot
\sum\limits_{k_0\leq |n|} |\varphi_n|+O(1),
\]
where $O(1)$ depends on $t_0$ and $\rho$. Since (\ref{phi}) this
series does not converge.

Thus condition $\varphi\in L_2^a(\mathbb{T}^N)$, $a>\frac{N}{2}$,
of Theorem \ref{main} cannot be replaced by the condition
(\ref{class}) in at least one-dimensional case.

In conclusion, note that a similar result with Theorem \ref{main}
for equation (\ref{eq}) with the Riemann-Liouville derivative is
valid for all functions $\varphi\in C(\mathbb{T}^N)\cap
L_2^{a-2}(\mathbb{T}^N)$ (see \cite{AO1}).

\

\bibliographystyle{amsplain}

\end{document}